\documentclass[oneside, reqno] {amsart}

\usepackage{xypic}
\input xy
\xyoption{all}
\usepackage{epsfig}
\usepackage{amsthm}
\usepackage{amssymb}
\usepackage{amsmath}
\usepackage{amscd}
\usepackage{color}
\usepackage[T1]{fontenc}
\usepackage{upgreek}
\usepackage[font=scriptsize]{caption}
\usepackage{wrapfig}
\usepackage{graphicx}
\usepackage{subfigure}

\newcommand{\bg}{\begin{equation}}
\newcommand{\ed}{\end{equation}}
\newcommand{\bga}{\begin{eqnarray}}
\newcommand{\eda}{\end{eqnarray}}
\newcommand{\pf}{\textbf{Proof:\ }}

\def\cbdu{\par{\raggedleft$\Box$\par}}

\newtheorem {Theorem}  {Theorem}

\numberwithin{Theorem}{section}

\newtheorem {Lemma}[Theorem]  {Lemma}
\newtheorem {Proposition}[Theorem]{Proposition}
\theoremstyle{definition}
\newtheorem{Definition}[Theorem]{Definition}
\theoremstyle{remark}

\expandafter\chardef\csname pre amssym.def
at\endcsname=\the\catcode`\@ \catcode`\@=11
\def\undefine#1{\let#1\undefined}
\def\newsymbol#1#2#3#4#5{\let\next@\relax
 \ifnum#2=\@ne\let\next@\msafam@\else
 \ifnum#2=\tw@\let\next@\msbfam@\fi\fi
 \mathchardef#1="#3\next@#4#5}
\def\mathhexbox@#1#2#3{\relax
 \ifmmode\mathpalette{}{\m@th\mathchar"#1#2#3}%
 \else\leavevmode\hbox{$\m@th\mathchar"#1#2#3$}\fi}
\def\hexnumber@#1{\ifcase#1 0\or 1\or 2\or 3\or 4\or 5\or 6\or 7\or 8\or
 9\or A\or B\or C\or D\or E\or F\fi}

\font\teneufm=eufm10 \font\seveneufm=eufm7 \font\fiveeufm=eufm5
\newfam\eufmfam
\textfont\eufmfam=\teneufm \scriptfont\eufmfam=\seveneufm
\scriptscriptfont\eufmfam=\fiveeufm

\catcode`\@=\csname pre amssym.def at\endcsname

\newcounter{remark}
\newcommand{\supp}{{\mathit supp}\,}

\def  \12  {{\frac{1}{2}}}

\def\build#1_#2^#3{\mathrel{\mathop{\kern 0pt#1}\limits_{#2}^{#3}}}

\numberwithin{equation}{section}

\begin{document}
%\currannalsline{0}{2006}

\title[Forced SQG]{Non-unique stationary solutions of forced SQG}

%\author{hello}

\author [Mimi Dai]{Mimi Dai}

\address{Department of Mathematics, Statistics and Computer Science, University of Illinois at Chicago, Chicago, IL 60607, USA}
\email{mdai@uic.edu}

\author [Qirui Peng]{Qirui Peng}

\address{Department of Mathematics, Statistics and Computer Science, University of Illinois at Chicago, Chicago, IL 60607, USA}
\email{qpeng9@uic.edu}

\thanks{The authors are partially supported by the NSF grant DMS--2009422. M. Dai is also supported by the AMS Centennial Fellowship.}

\begin{abstract}

We show the existence of non-unique stationary weak solutions for forced surface quasi-geostrophic (SQG) equation via a convex integration scheme. The scheme is implemented for the sum-difference system of two distinct solutions. Through this scheme, one observes the external forcing is naturally generated accompanying the flexibility in means of lack of uniqueness. It thus provides a transparent way to reveal the flexibility of the system with the presence of a forcing.

\bigskip


\hspace{0.02cm}CLASSIFICATION CODE: 35Q35, 35Q86, 76D03.
\end{abstract}

\maketitle

\section{Introduction}
\label{sec-int}
The two dimensional surface quasi-geostrophic equation (SQG) with external forcing
\begin{equation}\label{sqg}
\begin{split}
\partial_t\theta+ u \cdot\nabla \theta=&-\nu\Lambda^\gamma \theta+f, \\
u=\nabla^{\perp}\Lambda^{-1}\theta=&\left( -\mathcal R_2\theta, \mathcal R_1\theta\right)
\end{split}
\end{equation}
describes the evolution of the surface temperature $\theta$ in a rapidly rotating and stratified flow with velocity $u$ in the presence of buoyancy $f$.  Parameter $\nu\geq 0$ is the dissipation coefficient. The Zygmund operator $\Lambda$ is defined as $\Lambda=(-\Delta)^{\frac12}$; $R_1$ and $R_2$ are Riesz transforms. We assume $0<\gamma<\frac32$. System (\ref{sqg}) is posed on the spatial time domain $\mathbb T^2\times [0,\infty)$. It belongs to the family of active scalar equations with the non-local operator $T=: \nabla^{\perp}\Lambda^{-1}$ and drift velocity $u=T[\theta]$. Note the operator $T$ is odd in the sense that its Fourier symbol is odd; and $\nabla\cdot u= 0$.

Beside its significance in the study of atmosphere and oceanography,  the inviscid SQG (\ref{sqg}) with $\nu=0$ shares analogous features with the 3D Euler equation in various aspects. In the viscous case $\nu>0$, (\ref{sqg}) has the natural scaling: if $\theta(x,t)$ is a solution of (\ref{sqg}) with forcing $f(x,t)$, the rescaled temperature $\theta_\lambda(x,t)=\lambda^{\gamma-1}\theta(\lambda x,\lambda^\gamma t)$ is also a solution with rescaled forcing 
$f_\lambda(x,t)=\lambda^{2\gamma-1}f(\lambda x,\lambda^\gamma t)$. For appropriate forcing $f$ (for instance, $f=0$), system (\ref{sqg}) has the a priori estimate in $L^\infty(\mathbb T^2)$. While the space $L^\infty(\mathbb T^2)$ is scaling invariant for (\ref{sqg}) with $\gamma=1$, system (\ref{sqg}) is referred as critical with $\gamma=1$, supercritical for $\gamma<1$ and subcritical for $\gamma>1$. Since the early work \cite{Ped, CMT}, SQG has been extensively studied in the literature. Global existence of weak solutions with finite energy to the unforced SQG with $\nu\geq 0$ and $0<\gamma\leq 2$ was established by Resnick \cite{Res}. Global regular solution for the unforced critical SQG with $\gamma=1$ was obtained by the groups, Kieslev, Nazarov and Volberg \cite{KNV}, Caffarelli and Vasseur \cite{CaV} and Constantin and Vicol \cite{CV} applying different techniques. 

The forced SQG has also been investigated by many mathematicians. In particular, Kieslev and Nazarov \cite{KN} showed the existence of global regular solution to the critical SQG (\ref{sqg}) with the forcing of an ambient buoyancy gradient. For a special class of time independent forcing, Constantin, Tarfulea and Vicol \cite{CTV1, CTV2} studied the large time behavior of SQG solutions and proved the absence of anomalous dissipation. For external forcing $f\in L^p(\mathbb T^2)$ with $p>2$, Cheskidov and Dai \cite{CD} proved that the forced critical SQG has a compact global attractor in $L^2(\mathbb T^2)$. Regarding steady state of the forced critical SQG, Friedlander, Pavlovi\'c and Vicol \cite{FPV} showed that the steady state is nonlinearly unstable if the associated linear operator has spectrum in the unstable region. 

The main concern of this paper is on non-uniqueness of weak solutions for the forced SQG. In this line of research, active scalar equations with both even and odd drift operators have been investigated previously with the application of convex integration method, which was first ingeniously brought to Euler equation by De Lellis and Sz\'ekelyhidi \cite{DLS1, DLS2} from differential geometry. Shvydkoy \cite{Shv} showed the existence of non-unique bounded weak solutions for inviscid active scalar equations with even drift operator. Isett and Vicol \cite{IV} further studied the inviscid active scaler equation with a non-odd drift operator $T$ and obtained nontrivial compactly supported weak solutions with H\"older regularity $C^{\frac{1}{4d+1}-}$ on $\mathbb T^d$. For active scalar equations with odd operator $T$, including the SQG equation, the situation is different since the cancellation property due to the odd feature of $T$ presents a barrier to construct weak solutions with high regularity using convex integration scheme, see \cite{IV}. For the unforced SQG, Buckmaster, Shkoller and Vicol \cite{BSV} constructed nontrivial weak solutions with regularity $\Lambda^{-1}\theta\in C_t^\sigma C_x^\alpha$ for $\frac12<\alpha<\frac 45$ and $\sigma<\frac{\alpha}{2-\alpha}$, by working with the SQG momentum equation (the equation of $\Lambda^{-1} u$). Later on, by working directly with the $\theta$ equation, Isett and Ma \cite{IM} provided another construction of non-trivial solutions with the same H\"older regularity for the unforced SQG. On the other hand, working with the equation of the scalar function $\Lambda^{-1}\theta$, Cheng, Kwon and Li \cite{CKL} showed the existence of nontrivial stationary solutions to the unforced SQG with regularity $\Lambda^{-1}\theta\in C_x^\alpha$ for $\frac12<\alpha<\frac 23$ and $0<\gamma< 2-\alpha$. We observe that the constructed nontrivial stationary solution is less regular than the nontrivial time dependent solution. This is expected from the convex integration method since the temporal effect can play an important role to reduce certain errors in the iterative process.

\medskip

\subsection{Equivalent form of the forced system} 
\label{sec-two}
 
Our aim is to construct non-uniqune weak solutions to the forced stationary SQG, i.e. (\ref{sqg}) with $\partial_t\theta\equiv 0$, in space with higher H\"older regularity. The presence of the external forcing gives extra flexibility to the underlying problem. Such flexibility was revealed in the work \cite{DF} of the first author and Friedlander for forced dyadic models and in the paper \cite{FK} of Filonov and Khodunov as well. It is further elaborated in the following using the sum-difference reformulation of two distinct solutions appeared in \cite{FK}. It is easy to find an initial pair $(\theta, u, f_1)$ and $(\widetilde\theta, \widetilde u, f_2)$ satisfying (\ref{sqg}), with $\theta\neq \widetilde \theta$, $u=T[\theta]$ and $\widetilde u=T[\widetilde\theta]$ but without requiring $f_1\equiv f_2$. Denote the sum $p=\frac12(\theta+\widetilde\theta)$ and the difference $m=\frac12(\theta-\widetilde\theta)$, and hence $\theta=p+m$ and $\widetilde\theta=p-m$. The pair $(p, m)$ satisfies the forced system
\begin{equation}\label{pm2}
\begin{split}
p_t+T[p]\cdot\nabla p+T[m]\cdot\nabla m=&-\nu \Lambda^\gamma p+\frac12(f_1+f_2),\\
m_t+T[p]\cdot\nabla m+T[m]\cdot\nabla p=&-\nu\Lambda^\gamma m+\frac12(f_1-f_2),\\
\nabla\cdot T[p]= 0, \ \ \nabla\cdot T[m]=&\ 0.
\end{split}
\end{equation} 
The extra flexibility relies on the fact that forcing presents in both equations of (\ref{pm2}). If it were the case $f_1\equiv f_2=f$, we would have found two distinct solutions $(\theta, u, f)$ and $(\widetilde\theta, \widetilde u, f)$ of (\ref{sqg}) with the same forcing function. Naturally, to achieve the goal of obtaining two distinct solutions of the forced SQG, we treat $\frac12(f_1-f_2)$ as an initial error forcing term and apply a convex integration scheme to reduce it iteratively and eventually remove it. Three remarks are unfolded regarding the convex integration scheme in this context: (i) the convex integration scheme will be only applied to the $m$ equation not the $p$ equation in (\ref{pm2}); (ii) inspired by the work \cite{IM}, we recast the forcing $\frac12(f_1+f_2)=\Delta G$ and $\frac12(f_1-f_2)=\Delta \widetilde G$ into second derivative form; (iii) to improve the error estimates, we perform a special two-step construction in each iteration stage. Further details will be provided in Section \ref{sec-outline}.

\medskip

\subsection{Notion of weak solutions and main result}
\label{sec-result}

Since $\nabla\cdot u=0$, it is natural to define a stationary weak solution $\theta\in L^2(\mathbb T^2)$ of (\ref{sqg}) if the integral equation
\[-\int_{\mathbb T^2} \theta u\cdot \nabla\psi \, dx+ \nu\int_{\mathbb T^2} \theta \Lambda^{\gamma} \psi \, dx=\int_{\mathbb T^2} f \psi \, dx\]
holds for any $\psi\in C^{\infty}(\mathbb T^2)$. However, thanks to the odd feature of the operator $T$ for SQG, a weak solution to (\ref{sqg}) in the distributional sense can be defined in $\dot H^{-\frac12}$. Indeed, denote the commutator 
\[ [\nabla^{\perp} \Lambda^{-1}, \nabla\psi] f=\nabla^{\perp} \Lambda^{-1}(f\nabla\psi)-\nabla^{\perp} \Lambda^{-1}f\cdot \nabla\psi.\]
We observe that if $f\in \dot H^{-\frac12}$, $[\nabla^{\perp} \Lambda^{-1}, \nabla\psi] f\in \dot H^{\frac12}$.
On the other hand, integration by parts yields
\begin{equation}\notag
\begin{split}
\int_{\mathbb T^2} \theta u\cdot \nabla\psi \, dx=& \int_{\mathbb T^2} \theta \nabla^{\perp}\Lambda^{-1}\theta\cdot \nabla\psi \, dx\\
=& -\int_{\mathbb T^2} \theta \nabla^{\perp}\Lambda^{-1}\cdot(\theta \nabla\psi) \, dx\\
=&-\int_{\mathbb T^2} \theta u\cdot \nabla\psi \, dx-\int_{\mathbb T^2} \theta [\nabla^{\perp} \Lambda^{-1}, \nabla\psi]\theta\, dx.
\end{split}
\end{equation}
Thus we have
\begin{equation}\notag
\begin{split}
\int_{\mathbb T^2} \theta u\cdot \nabla\psi \, dx
=&-\frac12\int_{\mathbb T^2} \theta [\nabla^{\perp} \Lambda^{-1}, \nabla\psi]\theta\, dx
\end{split}
\end{equation}
and the right hand side integral is well-defined for $\theta\in\dot H^{-\frac12}$.

\begin{Definition}\label{def}
A distribution $\theta\in\dot H^{-\frac12}(\mathbb T^2)$ is said to be a stationary weak solution of (\ref{sqg}) with $f\in \dot H^{-r}$ if 
\begin{equation}\notag
\frac12\int_{\mathbb T^2}\Lambda^{-\frac12} \theta \Lambda^{\frac12}\left([\nabla^{\perp}\Lambda^{-1},\nabla\psi]\theta \right)\, dx+\nu \int_{\mathbb T^2}\Lambda^{-\frac12} \theta \Lambda^{\gamma+\frac12}\psi\, dx=\int_{\mathbb T^2}  \Lambda^{-r}f \Lambda^{r}\psi\, dx
\end{equation} 
holds for any smooth function $\psi\in C^\infty(\mathbb T^2)$.
\end{Definition}

Existence of weak solutions of (\ref{sqg}) in $\dot H^{-\frac12}$ without external forcing was established by Marchand \cite{Mar}. 

The main result of this paper is the non-uniqueness of stationary weak solutions to (\ref{sqg}) with certain external forcing.
\begin{Theorem}\label{thm}
Let $\nu\geq 0$, $0<\gamma<2-\alpha$ and $\frac12\leq \alpha<\frac34$. There exists $G\in C^{2\alpha-1}(\mathbb T^2)$ such that there are at least two stationary weak solutions $\theta, \widetilde \theta$ of (\ref{sqg}) with the external forcing $f=\Delta G$ and $\Lambda^{-1}\theta, \Lambda^{-1}\widetilde\theta\in C^{\alpha}(\mathbb T^2)$.
\end{Theorem}

We mention that stationary weak solutions constructed here for the forced SQG have higher regularity ($C^{-\frac14}$) than the solutions constructed in \cite{CKL} for the SQG not driven by any forcing. The latter ones have regularity $C^{-\frac13}$.

\medskip

We conclude the introduction by laying out the organization of the rest of the paper. In Section \ref{sec-outline}, we sketch a general convex integration scheme for forced equation in order to construct non-unique solutions. Section \ref{sec-proof} is devoted to the proof of Theorem \ref{thm}, where an iterative statement is the crucial element.

\bigskip

\section{A convex integration scheme for forced equation}
\label{sec-outline}

%Denote $\eta=\Lambda^{-1}\theta$. System (\ref{sqg}) is equivalent to 
%\begin{equation}\label{sqg2}
%\partial_t\Lambda\eta+ \nabla^{\perp}\eta \cdot\nabla \Lambda\eta=-\nu\Lambda^{\gamma+1}\eta+f.
%\end{equation}
We outline a scheme of iteration and approximation for the forced SQG in this section. The scheme is certainly generic and can be adapted to any forced equations regardless of being stationary or time-dependent. 

Adapting the notations 
\[\eta=\Lambda^{-1}\theta, \ \ \widetilde\eta=\Lambda^{-1}\widetilde\theta,\ \ \Pi=\frac12(\eta+\widetilde \eta), \ \ \mu=\frac12(\eta-\widetilde \eta),\] 
we have 
\[p=\Lambda \Pi, \ \ m=\Lambda \mu, \ \ T[p]=\nabla^{\perp}\Pi, \ \ T[m]=\nabla^{\perp}\mu.\]
Following the idea of \cite{CKL}, the convex integration is performed at the level of $\mu$. 
Thus (\ref{pm2}) can be written as
\begin{equation}\label{pm3}
\begin{split}
\partial_t \Lambda \Pi+\nabla\cdot(\Lambda \Pi\nabla^{\perp}\Pi)+\nabla\cdot(\Lambda \mu\nabla^{\perp}\mu)=&-\nu \Lambda^{\gamma+1} \Pi+\Delta G,\\
\partial_t \Lambda \mu+\nabla\cdot(\Lambda \mu\nabla^{\perp}\Pi)+\nabla\cdot(\Lambda \Pi\nabla^{\perp}\mu)=&-\nu\Lambda^{\gamma+1} \mu+\Delta \widetilde G\\
%\nabla\cdot \nabla^{\perp}\Pi= 0, \ \ \nabla\cdot \nabla^{\perp}\mu=&\ 0.
\end{split}
\end{equation} 
with forcing functions $G$ and $\widetilde G$ satisfying $\Delta G=\frac12(f_1+f_2)$ and $\Delta \widetilde G=\frac12(f_1-f_2)$.
We consider stationary solutions of (\ref{pm3}), i.e. solutions to the forced system
\begin{equation}\label{pm4}
\begin{split}
\Lambda \Pi\nabla^{\perp}\Pi+\Lambda \mu\nabla^{\perp}\mu=&-\nu \Lambda^{\gamma-1} \nabla\Pi+\nabla  G,\\
\Lambda \mu\nabla^{\perp}\Pi+\Lambda \Pi\nabla^{\perp}\mu=&-\nu\Lambda^{\gamma-1}\nabla \mu+\nabla\widetilde G.
%\nabla\cdot \nabla^{\perp}\Pi= 0, \ \ \nabla\cdot \nabla^{\perp}\mu=&\ 0.
\end{split}
\end{equation}

The goal is to construct a sequence of approximating solutions $\{(\Pi_q,\mu_q, G_q, \widetilde G_q)\}_{q\geq 0}$ of (\ref{pm4}) such that  $\widetilde G_q$ approaches zero in a suitable norm as $q\to\infty$. Thus the limit $(\Pi, \mu, G, 0)$ with non-vanishing $\mu$ is a solution of (\ref{pm4}). Equivalently, it implies the existence of two distinct stationary solutions $\theta=\Lambda(\Pi+\mu)$ and $\widetilde\theta=\Lambda(\Pi-\mu)$ of (\ref{sqg}) with forcing $f=\Delta G$.

We apply an iterative process to construct a sequence of approximating solutions. Let $q$ be even, starting from $q=0$. In general, let $(\Pi_q,\mu_q,  G_q, \widetilde G_q)$ be the solution of (\ref{pm4}) at the $q$-th iteration, i.e.
\begin{equation}\label{pm-q}
\begin{split}
\Lambda \Pi_q\nabla^{\perp}\Pi_q+\Lambda \mu_q\nabla^{\perp}\mu_q=&-\nu \Lambda^{\gamma-1} \nabla\Pi_q+\nabla G_q,\\
\Lambda \mu_q\nabla^{\perp}\Pi_q+\Lambda \Pi_q\nabla^{\perp}\mu_q=&-\nu\Lambda^{\gamma-1}\nabla \mu_q+\nabla \widetilde G_q
%\nabla\cdot \nabla^{\perp}\Pi_q= 0, \ \ \nabla\cdot \nabla^{\perp}\mu_q=&\ 0.
\end{split}
\end{equation} 
Recall 
\[\Pi_q=\frac12\Lambda^{-1}(\theta_q+\widetilde\theta_q), \ \ \mu_q=\frac12\Lambda^{-1}(\theta_q-\widetilde\theta_q)\]
and hence
\[\theta_q=\Lambda (\Pi_q+\mu_q), \ \ \widetilde\theta_q=\Lambda (\Pi_q-\mu_q), \ \ f_{1,q}=\Delta(G_q+\widetilde G_q), \ \ f_{2,q}=\Delta(G_q-\widetilde G_q). \]
Both $(\theta_q, f_{1,q})$ and $(\widetilde\theta_q, f_{2,q})$ satisfy the stationary forced SQG equation (\ref{sqg}). 
 Each stage of the construction consists two steps, from $q$-th to $(q+1)$-th step and from $(q+1)$-th to $(q+2)$-th step.  
 
In the first step, we construct $M_{q+1}$ to produce $(\Pi_{q+1}, \mu_{q+1})$
\[\mu_{q+1}=\mu_q+M_{q+1}, \ \ \Pi_{q+1}=\Pi_q-M_{q+1}.\]
We observe that 
\begin{equation}\notag
\begin{split}
\Lambda^{-1}\theta_{q+1}=&\ \eta_{q+1}= \Pi_{q+1}+\mu_{q+1}=\Pi_q+\mu_q\\
=&\ \eta_q=\Lambda^{-1}\theta_{q}, \\
\Lambda^{-1}\widetilde\theta_{q+1}=&\ \widetilde\eta_{q+1}= \Pi_{q+1}-\mu_{q+1}=\Pi_q-\mu_q-2M_{q+1}\\
=&\ \widetilde\eta_q-2M_{q+1}=\Lambda^{-1}\widetilde\theta_{q}-2M_{q+1}.
\end{split}
\end{equation}
Denote $G_{q+1}$ and $\widetilde G_{q+1}$ by the stress functions associated with $\Pi_{q+1}$ and $\mu_{q+1}$ respectively. The tuplet $(\Pi_{q+1}, \mu_{q+1}, G_{q+1}, \widetilde G_{q+1})$ satisfies 
\begin{equation}\label{pm-q1}
\begin{split}
\Lambda \Pi_{q+1}\nabla^{\perp}\Pi_{q+1}+\Lambda \mu_{q+1}\nabla^{\perp}\mu_{q+1}=&-\nu \Lambda^{\gamma-1} \nabla\Pi_{q+1}+\nabla G_{q+1},\\
\Lambda \mu_{q+1}\nabla^{\perp}\Pi_{q+1}+\Lambda \Pi_{q+1}\nabla^{\perp}\mu_{q+1}=&-\nu\Lambda^{\gamma-1}\nabla \mu_{q+1}+\nabla\widetilde G_{q+1}
%\nabla\cdot \nabla^{\perp}\Pi_{q+1}= 0, \ \ \nabla\cdot \nabla^{\perp}\mu_{q+1}=&\ 0.
\end{split}
\end{equation} 
Subtraction of the second equation of (\ref{pm-q}) from the second equation of (\ref{pm-q1}) leads to
\begin{equation}\label{g-q1}
\begin{split}
\nabla\widetilde G_{q+1}=&\ \nu\Lambda^{\gamma-1}\nabla M_{q+1}+\left(\Lambda \widetilde\eta_q\nabla^{\perp} M_{q+1}+ \Lambda M_{q+1}\nabla^{\perp} \widetilde\eta_q\right)\\
&+\left(\nabla\widetilde G_q-2\Lambda M_{q+1}\nabla^{\perp}M_{q+1}\right).
\end{split}
\end{equation}
To make a remark, $M_{q+1}$ will be constructed such that 
\begin{equation}\label{cancel1}
\nabla\widetilde G_q-2\Lambda M_{q+1}\nabla^{\perp}M_{q+1}\sim \nabla^{\perp} F_{q+1}
\end{equation}
with small error compared to other terms on the right hand side of (\ref{g-q1}) and some function $F_{q+1}$. 

In the second step of this stage, we construct $M_{q+2}$ such that 
\[\mu_{q+2}=\mu_{q+1}+M_{q+2}, \ \ \Pi_{q+2}=\Pi_{q+1}+M_{q+2}.\]
Again we have equivalently 
\begin{equation}\notag
\begin{split}
\Lambda^{-1}\theta_{q+2}=&\ \eta_{q+2}= \Pi_{q+2}+\mu_{q+2}=\Pi_{q+1}+\mu_{q+1}+2M_{q+2}\\
=&\ \eta_{q+1}+2M_{q+2}=\Lambda^{-1}\theta_{q+1}+2M_{q+2}, \\
\Lambda^{-1}\widetilde\theta_{q+2}=&\ \widetilde\eta_{q+2}= \Pi_{q+2}-\mu_{q+2}=\Pi_{q+1}-\mu_{q+1}\\
=&\ \widetilde\eta_{q+1}=\Lambda^{-1}\widetilde\theta_{q+1}.
\end{split}
\end{equation}
Let $\Pi_{q+2}$ and $\mu_{q+2}$ satisfy the system with functions $G_{q+2}$ and $\widetilde G_{q+2}$
\begin{equation}\label{pm-q2}
\begin{split}
\Lambda \Pi_{q+2}\nabla^{\perp}\Pi_{q+2}+\Lambda \mu_{q+2}\nabla^{\perp}\mu_{q+2}=&-\nu \Lambda^{\gamma-1} \nabla\Pi_{q+2}+\nabla G_{q+2},\\
\Lambda \mu_{q+2}\nabla^{\perp}\Pi_{q+2}+\Lambda \Pi_{q+2}\nabla^{\perp}\mu_{q+2}=&-\nu\Lambda^{\gamma-1}\nabla \mu_{q+2}+\nabla\widetilde G_{q+2}
%\nabla\cdot \nabla^{\perp}\Pi_{q+2}= 0, \ \ \nabla\cdot \nabla^{\perp}\mu_{q+2}=&\ 0.
\end{split}
\end{equation} 
Taking the subtraction of the second equation in (\ref{pm-q1}) and the second equation in (\ref{pm-q2}) gives
\begin{equation}\label{g-q2}
\begin{split}
\nabla\widetilde G_{q+2}=&\ \nu\Lambda^{\gamma-1}\nabla M_{q+2}+\left(\Lambda \eta_{q+1}\nabla^{\perp} M_{q+2}+ \Lambda M_{q+2}\nabla^{\perp} \eta_{q+1}\right)\\
&+\left(\nabla\widetilde G_{q+1}+2\Lambda M_{q+2}\nabla^{\perp}M_{q+2}\right).
\end{split}
\end{equation}
Analogously, $M_{q+2}$ will be constructed to reduce the size of $\widetilde G_{q+1}$ in the sense that 
\begin{equation}\label{cancel2}
\nabla\widetilde G_{q+1}+2\Lambda M_{q+2}\nabla^{\perp}M_{q+2}\sim \nabla^{\perp} F_{q+2}
\end{equation}
for some function $F_{q+2}$.
Iterating the stages described above for even $q\geq 2$, we obtain a sequence $\{(\Pi_{q+i}, \mu_{q+i}, G_{q+i}, \widetilde G_{q+i})\}_{i=1,2; q\geq 0}$ with 
$(\Pi_{q+i}, \mu_{q+i}, G_{q+i}, \widetilde G_{q+i})$ satisfying the forced systems (\ref{pm-q1}) and (\ref{pm-q2}). In particular, the force functions $\widetilde G_{q+1}$ and
$\widetilde G_{q+2}$ satisfy (\ref{g-q1}) and (\ref{g-q2}) respectively. 
We observe that in this iteration process,
\begin{equation}\label{skip}
\eta_{q+1}=\eta_q, \ \ \widetilde \eta_{q+2}=\widetilde \eta_{q+1}, \ \ \mbox{for any even} \ \ q\geq 0
\end{equation}
which is crucial to control the Nash errors in (\ref{g-q1}) and (\ref{g-q2}) and hence improve the regularity of the constructed weak solutions. 

On the other hand we notice the iteration gives the stress functions $G_{q+1}$ and $G_{q+2}$ in the sum equations
\begin{equation}\label{g1-q1}
\begin{split}
\nabla G_{q+1}=&-\nu\Lambda^{\gamma-1}\nabla M_{q+1}-\left(\Lambda \widetilde\eta_q\nabla^{\perp} M_{q+1}+ \Lambda M_{q+1}\nabla^{\perp} \widetilde\eta_q\right)\\
&+2\Lambda M_{q+1}\nabla^{\perp}M_{q+1}+\nabla G_q,
\end{split}
\end{equation}
\begin{equation}\label{g1-q2}
\begin{split}
\nabla G_{q+2}=&\ \nu\Lambda^{\gamma-1}\nabla M_{q+2}+\left(\Lambda \eta_{q+1}\nabla^{\perp} M_{q+2}+ \Lambda M_{q+2}\nabla^{\perp} \eta_{q+1}\right)\\
&+2\Lambda M_{q+2}\nabla^{\perp}M_{q+2}+\nabla G_{q+1}.
\end{split}
\end{equation}
Comparing (\ref{g-q1}) and (\ref{g1-q1}) we observe that the ``reduced'' amount of force from $\widetilde G_q$ to $\widetilde G_{q+1}$ is gained by $G_{q+1}$ from $G_{q}$. In contrast, we note from (\ref{g-q2}) and (\ref{g1-q2}) that both $\widetilde G_{q+1}$ and $G_{q+1}$ are ``reduced'' by the same amount of force to $\widetilde G_{q+2}$ and $G_{q+2}$ respectively.  

We mention that this convex integration scheme gives the same improvement for Nash error estimate as the alternating scheme used in \cite{BHP} for the forced Euler equation.

\bigskip

\section{Proof of the main theorem}
\label{sec-proof}

\medskip

\subsection{Nonlocal operators and auxiliary lemmas}
\label{sec-prep}

We collect some useful estimates for some nonlocal operators and a crucial algebraic lemma on the decomposition of a stress function. The proofs can be found in \cite{CKL}.
\begin{Lemma}\label{le-Leib}
Let $\lambda \xi\in \mathbb Z^2$ for $|\xi|=1$ and $g(x)=a(x)\cos(\lambda\xi\cdot x)$. We have
\begin{equation}\notag
\Lambda g=\lambda g+(\xi\cdot \nabla a)\sin(\lambda \xi\cdot x)+T_{1, \lambda\xi}[a]\cos(\lambda \xi\cdot x)+T_{2, \lambda\xi}[a]\sin(\lambda \xi\cdot x)
\end{equation}
with
\begin{equation}\notag
\begin{split}
\widehat{T_{1, \lambda\xi}[a]}(k)=&\left(\frac12\left(|\lambda\xi+k|+|\lambda\xi-k|\right)-\lambda \right)\hat a(k),\\
\widehat{T_{2, \lambda\xi}[a]}(k)=&\ i\left(\frac12\left(|\lambda\xi+k|-|\lambda\xi-k|\right)-\xi\cdot k \right)\hat a(k).
\end{split}
\end{equation}
\end{Lemma}

\begin{Lemma}\label{le-aux1}
Let $a\in L^\infty(\mathbb T^2)$ with zero mean and $\supp(\hat a)\subset \{|k|\leq r\}$ for $\mu\geq 10$. Let $T$ be the Fourier multiplier defined by $\widehat {T[f]}(k)=m(k)\hat f(k)$ for a homogeneous function $m\in C^\infty(\mathbb R^2/ \{0\})$ of degree 0. Then 
\[\|T[a]\|_{L^\infty}\lesssim \|a\|_{L^\infty} \log r\]
up to a constant depending on $a$. 
\end{Lemma}

The notation $A\lesssim B$ represents an estimate up to a constant, that is, $A\leq C B$ for some $C>0$. It will be used often throughout the text when the implicit constant does not play a role.
 
%Let $f\in C^\infty(\mathbb T^2)$ and $m(\xi)=\int_{\mathbb R^2} K(z)e^{-i\xi\cdot z}\, dz$ for $K\in L^1(\mathbb R^2)$. Define 
%\begin{equation}\notag
%T[f](x)=\sum_{k}m(k)\hat f(k)e^{ik\cdot x}=\int_{\mathbb R^2} K(z)f(x-z)\, dz.
%\end{equation}
%Then we have
%\begin{equation}\notag
%\|T[f]\|_{L^p(\mathbb T^2)}\leq \|K\|_{L^1(\mathbb R^2)}\|f\|_{L^p(\mathbb T^2)}, \ \ 1\leq p\leq \infty.
%\end{equation}

%Let $f,g\in C^\infty(\mathbb T^2)$ and $m(\xi, \eta)=\int_{\mathbb R^2}\int_{\mathbb R^2} K(z_1,z_2)e^{-i\xi\cdot z_1-i\eta\cdot z_2}\, dz_1dz_2$ for $K\in L^1(\mathbb R^2\times \mathbb R^2)$. Define 
%\begin{equation}\notag
%\begin{split}
%T[f, g](x)=&\ \sum_{k}\left( \sum_{k'\in\mathbb Z^2}m(k', k-k')\hat f(k')\hat g(k-k') \right)e^{ik\cdot x}\\
%=&\int_{\mathbb R^2}\int_{\mathbb R^2} K(z_1, z_2)f(x-z_1)g(x-z_2)\, dz_1dz_2.
%\end{split}
%\end{equation}
%Then we have
%\begin{equation}\notag
%\|T[f,g]\|_{L^r(\mathbb T^2)}\leq \|K\|_{L^1(\mathbb R^2\times \mathbb R^2)}\|f\|_{L^p(\mathbb T^2)}\|g\|_{L^q(\mathbb T^2)}, \ \ 1\leq r, p, q\leq \infty
%\end{equation}
%with $\frac1r=\frac1p+\frac1q$.

\begin{Lemma}\label{le-t1}
Let $a_0: \mathbb T^2\to \mathbb R$ with $\supp(\hat a_0)\subset \{|k|\leq r\}$ for $10\leq r\leq \frac12\lambda$. We have
\begin{equation}\notag
\begin{split}
\|T_{1,\lambda\xi} [a_0]\|_{L^\infty} \lesssim&\ \lambda^{-1}r^2\|a_0\|_{L^\infty},\\
\|T_{2,\lambda\xi} [a_0]\|_{L^\infty} \lesssim&\ \lambda^{-2}r^3\|a_0\|_{L^\infty},\\
\|\Delta^{-1}\nabla T_{2,\lambda\xi} [a_0]\|_{X} \lesssim&\ \lambda^{-2}r^2\|a_0\|_{L^\infty} \log r.
\end{split}
\end{equation}
\end{Lemma}

\begin{Lemma}\label{le-t2}
Let $a_0$ be as in Lemma \ref{le-t1}. Let $m_{j, \ell}=m_{j, \ell, \lambda, \mu}$ with $j=1,2,3$ be the Fourier symbols defined by 
\begin{equation}\notag
\begin{split}
a_0T_{1, \lambda\xi}[a_0]=&\ \frac{r^2}{\lambda^2}\sum_{\ell=1}^2\partial_{x_{\ell}} T_{m_{1,\ell}}[a_0,a_0],\\
T_{1,\lambda\xi}[a_0] \partial_{x_1} a_0=&\ \frac{r^2}{\lambda}\sum_{\ell=1}^2\partial_{x_{\ell}}T_{m_{2,\ell}}[a_0,a_0],\\
T_{1,\lambda\xi}[a_0]\partial_{x_2}a_0=&\ \frac{r^2}{\lambda}\sum_{\ell=1}^2\partial_{x_\ell}T_{m_{3,\ell}}[a_0,a_0],
\end{split}
\end{equation}
and $K_{j,\ell}=\mathcal F^{-1}(m_{j,\ell})$. Then we have
\[\|K_{j,\ell}\|_{L^1(\mathbb R^4)}\lesssim 1\]
up to a constant independent of $\lambda$ and $r$.
\end{Lemma}

\begin{Lemma}[Algebraic Lemma]\label{le-alg}
Let $\xi_1=\left(\frac35, \frac45\right)$ and $\xi_2=(1,0)$. The Riesz transforms $\mathcal R_j^o$ with $j=1,2$ have the Fourier symbols 
\begin{equation}\notag
\widehat {\mathcal R_1^o}(k_1, k_2)=\frac{25(k_2^2-k_1^2)}{12|k|^2}, \ \ \widehat {\mathcal R_2^o}(k_1, k_2)=\frac{7(k_2^2-k_1^2)}{12|k|^2}+\frac{4k_1k_2}{|k|^2}.
\end{equation}
Then for any function $G\in C^\infty_0(\mathbb T^2)$, the decomposition 
\begin{equation}\notag
\nabla G=\sum_{j=1}^2\xi_j^{\perp}(\xi_j\cdot\nabla) (\mathcal R_j^o G)+\nabla^{\perp} F
\end{equation}
holds for some $F$.
\end{Lemma}

\medskip

\subsection{Building blocks}
\label{sec-blocks}

We consider the increment $M_{n+1}$ in the form
\begin{equation}\label{Mn}
M_{n+1}(x)=\sum_{j=1}^2a_{j,n+1}(x)\cos (\lambda_{n+1}\xi_j\cdot x)
\end{equation}
where $a_{j,n+1}(x)=:a_j(x)$ are magnitude functions to be determined in the following. 

In view of Lemma \ref{le-Leib}, it follows immediately that
\begin{equation}\notag
\begin{split}
\Lambda M_{n+1}=&\ \lambda_{n+1} M_{n+1}+\sum_{j=1}^2(\xi_j\cdot \nabla a_j)\sin(\lambda_{n+1}\xi_j\cdot x)\\
&+\sum_{j=1}^2T_{1, \lambda_{n+1}\xi_j}[a_j]\cos(\lambda_{n+1}\xi_j\cdot x)+\sum_{j=1}^2T_{2, \lambda_{n+1}\xi_j}[a_j]\sin(\lambda_{n+1}\xi_j\cdot x).
\end{split}
\end{equation}
On the other hand, we have
\begin{equation}\notag
\nabla^{\perp} M_{n+1}=\sum_{j=1}^2\nabla^{\perp} a_j\cos(\lambda_{n+1}\xi_j\cdot x)-\lambda_{n+1} \sum_{j=1}^2\xi_j^{\perp} a_j\sin(\lambda_{n+1}\xi_j\cdot x).
\end{equation}
Hence straightforward computation shows that
\begin{equation}\notag
\begin{split}
&\Lambda M_{n+1}\nabla^{\perp} M_{n+1}\\
=&\ \lambda_{n+1} M_{n+1}\nabla^{\perp} M_{n+1}-\lambda_{n+1} \sum_{1\leq j,j'\leq 2}(\xi_j\cdot\nabla a_{j})(\xi_{j'})^{\perp} a_{j'}\sin(\lambda_{n+1}\xi_j\cdot x)\sin(\lambda_{n+1}\xi_{j'}\cdot x)\\
&+\sum_{1\leq j,j'\leq 2}(\xi_j\cdot\nabla a_{j})\nabla^{\perp} a_{j'}\sin(\lambda_{n+1}\xi_j\cdot x)\cos(\lambda_{n+1}\xi_{j'}\cdot x)\\
&-\lambda_{n+1} \sum_{1\leq j,j'\leq 2}T_{1, \lambda_{n+1}\xi_j}[a_j](\xi_{j'})^{\perp} a_{j'}\cos(\lambda_{n+1}\xi_j\cdot x)\sin(\lambda_{n+1}\xi_{j'}\cdot x)\\
&-\lambda_{n+1} \sum_{1\leq j,j'\leq 2}T_{2, \lambda_{n+1}\xi_j}[a_j](\xi_{j'})^{\perp} a_{j'}\sin(\lambda_{n+1}\xi_j\cdot x)\sin(\lambda_{n+1}\xi_{j'}\cdot x)\\
&+ \sum_{1\leq j,j'\leq 2}T_{1, \lambda_{n+1}\xi_j}[a_j]\nabla^{\perp} a_{j'}\cos(\lambda_{n+1}\xi_j\cdot x)\cos(\lambda_{n+1}\xi_{j'}\cdot x)\\
&+ \sum_{1\leq j,j'\leq 2}T_{2, \lambda_{n+1}\xi_j}[a_j]\nabla^{\perp} a_{j'}\sin(\lambda_{n+1}\xi_j\cdot x)\cos(\lambda_{n+1}\xi_{j'}\cdot x)\\
=:&\ \frac12\lambda_{n+1} \nabla^{\perp} M_{n+1}^2+J_1+J_2+J_3+J_4+J_5+J_6.
\end{split}
\end{equation}
We further analyze the terms $J_1$, $J_4$ and $J_5$ by separating self interactions of plane waves from non-self interactions,
\begin{equation}\notag
\begin{split}
J_1=&-\lambda_{n+1} \sum_{j=1}^2(\xi_j\cdot\nabla a_{j})(\xi_j)^{\perp} a_{j}\sin^2(\lambda_{n+1}\xi_j\cdot x)\\
&-\lambda_{n+1} \sum_{j\neq j'}(\xi_j\cdot\nabla a_{j})(\xi_{j'})^{\perp} a_{j'}\sin(\lambda_{n+1}\xi_j\cdot x)\sin(\lambda_{n+1}\xi_{j'}\cdot x)\\
=&-\frac14\lambda_{n+1}\sum_{j=1}^2\xi_j^{\perp}(\xi_j\cdot\nabla)a^2_{j}+\frac12\lambda_{n+1} \sum_{j=1}^2\xi_j^{\perp}(\xi_j\cdot\nabla)a^2_{j}\cos(2\lambda_{n+1}\xi_j\cdot x)\\
&-\lambda_{n+1} \sum_{j\neq j'}(\xi_j\cdot\nabla a_{j})(\xi_{j'})^{\perp} a_{j'}\sin(\lambda_{n+1}\xi_j\cdot x)\sin(\lambda_{n+1}\xi_{j'}\cdot x),
\end{split}
\end{equation}
\begin{equation}\notag
\begin{split}
J_4=&-\lambda_{n+1}\sum_{j=1}^2T_{2, \lambda_{n+1}\xi_j}[a_j]\xi_j^{\perp} a_{j}\sin^2(\lambda_{n+1}\xi_j\cdot x)\\
&-\lambda_{n+1} \sum_{j\neq j'}T_{2, \lambda_{n+1}\xi_j}[a_j](\xi_{j'})^{\perp} a_{j'}\sin(\lambda_{n+1}\xi_j\cdot x)\sin(\lambda_{n+1}\xi_{j'}\cdot x)\\
=&-\frac12\lambda_{n+1}\sum_{j=1}^2T_{2, \lambda_{n+1}\xi_j}[a_j]\xi_j^{\perp} a_{j}+\frac12\lambda_{n+1}\sum_{j=1}^2T_{2, \lambda_{n+1}\xi_j}[a_j]\xi_j^{\perp} a_{j}\cos(2\lambda_{n+1}\xi_j\cdot x)\\
&-\lambda_{n+1} \sum_{j\neq j'}T_{2, \lambda_{n+1}\xi_j}[a_j](\xi_{j'})^{\perp} a_{j'}\sin(\lambda_{n+1}\xi_j\cdot x)\sin(\lambda_{n+1}\xi_{j'}\cdot x),
\end{split}
\end{equation}
\begin{equation}\notag
\begin{split}
J_5=&\ \frac12\sum_{j=1}^2T_{1, \lambda_{n+1}\xi_j}[a_j]\nabla^{\perp} a_{\xi_j}+\frac12\sum_{j=1}^2T_{1, \lambda_{n+1}\xi_j}[a_j]\nabla^{\perp} a_{j}\cos(2\lambda_{n+1}\xi_j\cdot x)\\
&+\sum_{j\neq j'}T_{1, \lambda_{n+1}\xi_j}[a_j]\nabla^{\perp} a_{j'}\cos(\lambda_{n+1}\xi_j\cdot x)\cos(\lambda_{n+1}\xi_{j'}\cdot x)\\
\end{split}
\end{equation}
where we note non-oscillatory terms are generated. Separation of self interactions from non-self interactions in $J_2$, $J_3$ and $J_6$ leads to
\begin{equation}\notag
\begin{split}
J_2=&\ \frac12\sum_{j=1}^2(\xi_j\cdot\nabla a_{j})\nabla^{\perp} a_{j}\sin(2\lambda_{n+1}\xi_j\cdot x)\\
&+\sum_{j\neq j'}(\xi_j\cdot\nabla a_{j})\nabla^{\perp} a_{j'}\sin(\lambda_{n+1}\xi_j\cdot x)\cos(\lambda_{n+1}\xi_{j'}\cdot x),
\end{split}
\end{equation}
\begin{equation}\notag
\begin{split}
J_3=& -\frac12\lambda_{n+1} \sum_{j=1}^2T_{1, \lambda_{n+1}\xi_j}[a_j]\xi_j^{\perp} a_{j}\sin(2\lambda_{n+1}\xi_j\cdot x)\\
&-\lambda_{n+1} \sum_{j\neq j'}T_{1, \lambda_{n+1}\xi_j}[a_j](\xi_{j'})^{\perp} a_{j'}\cos(\lambda_{n+1}\xi_j\cdot x)\sin(\lambda_{n+1}\xi_{j'}\cdot x),
\end{split}
\end{equation}
\begin{equation}\notag
\begin{split}
J_6=&\ \frac12\sum_{j=1}^2T_{2, \lambda_{n+1}\xi_j}[a_j]\nabla^{\perp} a_{j}\sin(2\lambda_{n+1}\xi_j\cdot x)\\
&+\sum_{j\neq j'}T_{2, \lambda_{n+1}\xi_j}[a_j]\nabla^{\perp} a_{j'}\sin(\lambda_{n+1}\xi_j\cdot x)\cos(\lambda_{n+1}\xi_{j'}\cdot x).
\end{split}
\end{equation}
Combining the algebra manipulations above we obtain
\begin{equation}\notag
\begin{split}
&\Lambda M_{n+1}\nabla^{\perp} M_{n+1}\\
=&\ \frac12\lambda_{n+1} \nabla^{\perp} M_{n+1}^2-\frac14\lambda_{n+1} \sum_{j=1}^2\xi_j^{\perp}(\xi_j\cdot\nabla)a^2_{j}\\
&-\frac12\lambda_{n+1}\sum_{j=1}^2T_{2, \lambda_{n+1}\xi}[a_j]\xi_j^{\perp} a_{\xi_j}+\frac12\sum_{j=1}^2T_{1, \lambda_{n+1}\xi_j}[a_j]\nabla^{\perp} a_{j}\\
&+\frac12\lambda_{n+1} \sum_{j=1}^2\left( \xi_j^{\perp}(\xi_j\cdot\nabla)a^2_{j}+T_{2, \lambda_{n+1}\xi_j}[a_j]\xi_j^{\perp} a_{j}+\lambda_{n+1}^{-1} T_{1, \lambda_{n+1}\xi_j}[a_j]\nabla^{\perp} a_{j} \right)\cos(2\lambda_{n+1}\xi_j\cdot x)\\
&+\frac12 \sum_{j=1}^2\left( (\xi_j\cdot\nabla a_{j})\nabla^{\perp} a_{j}-\lambda_{n+1}T_{1, \lambda_{n+1}\xi_j}[a_j]\xi_j^{\perp} a_{j}+T_{2, \lambda_{n+1}\xi_j}[a_j]\nabla^{\perp} a_{j} \right)\sin(2\lambda_{n+1}\xi_j\cdot x)\\
&-\lambda_{n+1} \sum_{j\neq j'}\left(\xi_j\cdot\nabla a_{j}+ T_{2, \lambda_{n+1}\xi_j}[a_j]\right)(\xi_{j'})^{\perp} a_{j'}\sin(\lambda_{n+1}\xi_j\cdot x)\sin(\lambda_{n+1}\xi_{j'}\cdot x)\\
&+\sum_{j\neq j'}\left(\xi_j\cdot\nabla a_{j}+T_{2, \lambda_{n+1}\xi_j}[a_j]\right)\nabla^{\perp} a_{j'}\sin(\lambda_{n+1}\xi_j\cdot x)\cos(\lambda_{n+1}\xi_{j'}\cdot x)\\
&-\lambda_{n+1} \sum_{j\neq j'}T_{1, \lambda_{n+1}\xi_j}[a_j](\xi_{j'})^{\perp} a_{j'}\cos(\lambda_{n+1}\xi_j\cdot x)\sin(\lambda_{n+1}\xi_{j'}\cdot x)\\
&+\sum_{j\neq j'}T_{1, \lambda_{n+1}\xi_j}[a_j]\nabla^{\perp} a_{j'}\cos(\lambda_{n+1}\xi_j\cdot x)\cos(\lambda_{n+1}\xi_{j'}\cdot x)\\
=:&\ \frac12\lambda_{n+1} \nabla^{\perp} M_{n+1}^2-\frac14\lambda_{n+1} \sum_{j=1}^2\xi_j^{\perp}(\xi_j\cdot\nabla)a^2_{j}\\
&+\nabla J_{NO}+\nabla J_{O1}+\nabla J_{O2}+\nabla J_{O3}+\nabla J_{O4}+\nabla J_{O5}+\nabla J_{O6}
\end{split}
\end{equation}
where $J_{NO}$ refers the non-oscillatory error in the second line and $J_{O1}, ... , J_{O6}$ the oscillatory errors in the order of the lines. Invoking Lemma \ref{le-alg},  the magnitude functions $a_j$ will be designed such that the major term 
\[-\frac14\lambda_{n+1} \sum_{j=1}^2\xi_j^{\perp}(\xi_j\cdot\nabla)a^2_{j}\]
cancels the principal part of the stress function in the sense of (\ref{cancel1}) and (\ref{cancel2}). Thus at a rough level, we expect to choose $a_{j}$ such that
\[\lambda_{n+1} a_{j}^2\sim \mathcal R_j^o \widetilde G \ \ \Longrightarrow \ \ a_{j}\sim \lambda_{n+1}^{-\frac12} (\mathcal R_j^o \widetilde G)^{\frac12}.\]
Before giving a precise definition of $a_{j}$, we carry out another heuristic argument on the sizes of $M_{n+1}$ and $\widetilde G_{n}$ in (\ref{cancel1}) and (\ref{cancel2}) with $n=q$ and $n=q+1$ respectively.  Fix a large constant $\lambda_0>0$. Let $b>1$. Choose the frequency number as the integer
\[\lambda_n=\left\lceil \lambda_0^{b^n}\right\rceil, \ \ \ n\in \mathbb N\cup \{0\}.\]
The magnitude measure is given by $\delta_n=\lambda_n^{-\beta}$ for a parameter $\beta>0$ to be specified later. We also choose the frequency localization number $r_{n+1}=(\lambda_n\lambda_{n+1})^{\frac12}$. In the process of iteration described in Subsection \ref{sec-outline}, we expect to have 
\[|\widetilde G_n|\sim \delta_n.\]
The increment $M_{n+1}$ will be constructed such that (i) it is supported in Fourier space near the frequency $\lambda_{n+1}$; (ii) it has $C^{\alpha}$ regularity for some $\alpha$ to be determined. 
Thus in view of the cancelations in (\ref{cancel1}) and (\ref{cancel2}), we have
\begin{equation}\label{m-q1}
|M_{n+1}|\sim \left(\lambda_{n+1}^{-1}\delta_{n}\right)^{\frac12}.
\end{equation}
The $C^{\alpha}$ regularity requirement for $M_{n+1}$ indicates 
\[\lambda_{n+1}^{\alpha} \left(\lambda_{n+1}^{-1}\delta_{n}\right)^{\frac12}\lesssim 1\]
which implies 
\begin{equation}\label{alpha}
\alpha<\frac12+\frac{\beta}{2b}.
\end{equation}
To realize the cancellations (\ref{cancel1}) and (\ref{cancel2}) we choose 
\begin{equation}\label{a}
\begin{split}
a_{j, q+1}=&\ \left(\frac{2\delta_q}{5\lambda_{q+1}}\right)^{\frac12}\left(c_0-\mathcal R_j^o\left(\frac{\widetilde G_q}{\delta_q}\right) \right)^{\frac12},\\
a_{j, q+2}=&\ \left(\frac{2\delta_{q+1}}{5\lambda_{q+2}}\right)^{\frac12}\left(c_0+\mathcal R_j^o\left(\frac{\widetilde G_{q+1}}{\delta_{q+1}}\right) \right)^{\frac12}
\end{split}
\end{equation}
where $c_0\geq 2$ is a constant such that the quantity in $(\cdot)^{\frac12}$ is positive. 
%For brevity, we denote $a_{\xi_j, q+1}=a_{j, q+1}$ and $a_{\xi_j, q+2}=a_{j, q+2}$. 
We then construct the increments $M_{q+1}$ and $M_{q+2}$ as 
\begin{equation}\label{m-q1-construct}
\begin{split}
M_{q+1}(x)=& \sum_{j=1}^2P_{\leq r_{q+1}} \left(a_{j, q+1}(x)\right)\cos(5\lambda_{q+1}\xi_j\cdot x),\\
M_{q+2}(x)=& \sum_{j=1}^2P_{\leq r_{q+2}} \left(a_{j, q+2}(x)\right)\cos(5\lambda_{q+2}\xi_j\cdot x)
\end{split}
\end{equation}
where $P_{\leq\mu_{n+1}}$ is a standard Littlewood-Paley projection operator.

\medskip

\subsection{Main iteration result}
Denote $X$ by the space of functions with the norm 
\[\|G\|_{X}=\|G\|_{L^\infty}+\sum_{j=1}^2\|\mathcal R_j^o G\|_{L^\infty}.\]

\begin{Proposition}\label{prop}
Let the parameters satisfy
\[ \lambda_0\gg1, \ \ \nu\geq 0, \ \ 0<\gamma<2-\alpha, \ \ \frac12\leq \alpha<\frac34.\]
There exist $b>1$ and $0<\beta<\frac12$ such that 
\begin{equation}\label{parameter}
(2\alpha-1)b<\beta< \min\left\{ \frac{2b}{2b-1}(\frac32-\gamma), \frac{b^2-2+2\alpha}{b(2b-1)} \right\}
\end{equation}
and the following holds. If $(\Pi_n, \mu_n, G_n, \widetilde G_n)$ satisfies (\ref{pm4}) with $\Pi_n, \mu_n\in C^\alpha$ and
%$\Pi_n$ and $\mu_n$ are localized to $\leq 6\lambda_n$, $G_n$ and $\widetilde G_n$ are localized to $\leq 12\lambda_n$, 
\begin{equation}\label{iter-1}
\Pi_n=P_{\leq 6\lambda_n}\Pi_n, \ \ \mu_n=P_{\leq 6\lambda_n}\mu_n, \ \ G_n=P_{\leq 12\lambda_n}G_n, \ \ \widetilde G_n=P_{\leq 12\lambda_n}\widetilde G_n,
\end{equation}
\begin{equation}\label{iter-2}
\|G_n\|_{X}\leq 1-\delta_n^{\frac12},
\end{equation}
\begin{equation}\label{iter-3}
\|G_n\|_{C^{s}}\lesssim \lambda_n^s \delta_n, \ \ \ s\geq \beta,
\end{equation}
\begin{equation}\label{iter-4}
\|\widetilde G_n\|_{X}\leq \delta_n.
\end{equation}
Then there exits $(\Pi_{n+1}, \mu_{n+1}, G_{n+1}, \widetilde G_{n+1})$ satisfying (\ref{pm4}) with $\Pi_{n+1}, \mu_{n+1}\in C^\alpha$, and (\ref{iter-1})-(\ref{iter-4}) satisfied with $n$ replaced by $n+1$.
%\[\|\Pi_{n+1}\|_{C^\alpha(\mathbb T^2)}+\|\mu_{n+1}\|_{C^\alpha(\mathbb T^2)}\leq 200, \ \ \|\widetilde G_{n+1}\|_{X}\leq \delta_{n+1},\ \ \|G_{n+1}\|_{C^{2\alpha-1}}\leq 100.\]
\end{Proposition}

\pf
First of all, one can check there exsit $b>1$ and $0<\beta<\frac12$ such that the parameter conditions in (\ref{parameter}) are satisfied. 
As in the iterative scheme sketched in Subsection \ref{sec-outline}, we need to prove the statement for $n=q$ and $n=q+1$ for any even integer $q\geq 0$. When $n=q$, we construct $M_{n+1}=M_{q+1}$ as appeared in (\ref{m-q1-construct}). Let 
\[\mu_{n+1}=\mu_n+M_{n+1}, \ \ \Pi_{n+1}=\Pi_n-M_{n+1}.\]
As pointed out in Subsection \ref{sec-outline} we have 
\begin{equation}\notag
\begin{split}
\Lambda^{-1}\theta_{n+1}=&\ \eta_{n+1}= \eta_n=\Lambda^{-1}\theta_{n}, \\
\Lambda^{-1}\widetilde\theta_{n+1}=&\ \widetilde\eta_{n+1}=\widetilde\eta_n-2M_{n+1}=\Lambda^{-1}\widetilde\theta_{n}-2M_{n+1}.
\end{split}
\end{equation}
For $G_{n+1}$ and $\widetilde G_{n+1}$ defined as in (\ref{g1-q1}) and (\ref{g-q1}) respectively with $q=n$, the tuplet $(\Pi_{n+1}, \mu_{n+1}, G_{n+1}, \widetilde G_{n+1})$ satisfies (\ref{pm4}), and (\ref{iter-1}) holds with $n$ replaced by $n+1$. In view of (\ref{a}) and (\ref{m-q1-construct}), we have
\begin{equation}\notag
\|M_{n+1}\|_{C^\alpha}\lesssim \lambda_{n+1}^\alpha \delta_n^{\frac12}\lambda_{n+1}^{-\frac12}\lesssim \lambda_{n+1}^{\alpha-\frac{\beta}{2b}-\frac12}\lesssim 1
\end{equation}
since $(2\alpha-1)b<\beta$. It follows immediately from $\Pi_n, \mu_n\in C^\alpha$ that $\Pi_{n+1}, \mu_{n+1}\in C^\alpha$.

We are left to show the estimates for $G_{n+1}$ and $\widetilde G_{n+1}$. We recall (\ref{g-q1})
\begin{equation}\label{est-gn1}
\begin{split}
\nabla\widetilde G_{n+1}=&\ \nu\Lambda^{\gamma-1}\nabla M_{n+1}+\left(\Lambda \widetilde\eta_n\nabla^{\perp} M_{n+1}+ \Lambda M_{n+1}\nabla^{\perp} \widetilde\eta_n\right)\\
&+\left(\nabla\widetilde G_n-2\Lambda M_{n+1}\nabla^{\perp}M_{n+1}\right)\\
=:&\ \nabla \widetilde G_D+\nabla \widetilde G_N+\nabla \widetilde G_R
\end{split}
\end{equation}
with $\widetilde G_D$, $\widetilde G_N$ and $\widetilde G_R$ denoting the dissipation error, Nash error and reduced error accordingly. We estimate the errors in the following. 

We can choose $\widetilde G_D=\nu\Lambda^{\gamma-1} M_{n+1}\in C^\infty_0(\mathbb T^2)$.
It follows from (\ref{a})
\begin{equation}\notag
\|a_{j,n+1}\|_{L^\infty}\lesssim \left( \frac{\delta_n}{\lambda_{n+1}}\right)^{\frac12},
\end{equation}
and hence 
\begin{equation}\notag
\|M_{n+1}\|_{L^\infty}\leq \sum_{j=1}^2\|a_{j,n+1}\|_{L^\infty} \lesssim \left( \frac{\delta_n}{\lambda_{n+1}}\right)^{\frac12}.
\end{equation}
Thus
\begin{equation}\label{est-gn2}
\|\widetilde G_D\|_{X}\lesssim \lambda_{n+1}^{\gamma-1}\|M_{n+1}\|_{L^\infty} \lesssim \lambda_{n+1}^{\gamma-1}\left( \frac{\delta_n}{\lambda_{n+1}}\right)^{\frac12}\sim\lambda_{n+1}^{\gamma-\frac32-\frac{\beta}{2b}}\leq \frac13\delta_{n+1}
\end{equation}
%provided $\gamma-\frac32-\frac{\beta}{2b}\leq -\beta$ which is equivalent to $\gamma<\frac32-\beta+\frac{\beta}{2b}$.
where the last step holds thanks to $\lambda_0\gg1$ and $\beta<\frac{2b}{2b-1}(\frac32-\gamma)$.

Regarding the Nash error we choose 
\[\widetilde G_N= \Delta^{-1}\nabla\cdot \left(\Lambda \widetilde\eta_n\nabla^{\perp} M_{n+1}+ \Lambda M_{n+1}\nabla^{\perp} \widetilde\eta_n\right).\]
In view of (\ref{skip}), $\widetilde\eta_n=\widetilde\eta_{n-1}=\Pi_{n-1}-\mu_{n-1}$. Since $\Pi_{n-1}, \mu_{n-1}\in C^\alpha$ based on iteration,
%\[\|\Pi_{n-1}\|_{C^\alpha}+\|\mu_{n-1}\|_{C^\alpha}\leq 200,\]
we have
\[\|\widetilde\eta_{n-1}\|_{L^\infty}\lesssim \|\Pi_{n-1}\|_{L^\infty}+\|\mu_{n-1}\|_{L^\infty}\lesssim \lambda_{n-1}^{-\alpha}.\] 
Therefore, it follows
\begin{equation}\label{est-gn3}
\begin{split}
\|\widetilde G_N\|_{X}\lesssim &\ \|M_{n+1}\|_{L^\infty}\left(\|\nabla^{\perp}\widetilde \eta_{n-1}\|_{L^\infty}+\|\Lambda \widetilde\eta_{n-1}\|_{L^\infty} \right)\\
\lesssim& \ \delta_n^{\frac12}\lambda_{n+1}^{-\frac12}\lambda_{n-1}^{1-\alpha}\\
\leq&\ \frac13 \delta_{n+1}
\end{split}
\end{equation}
provided 
\begin{equation}\notag
-\frac12\beta-\frac12b+\frac{1}{b}(1-\alpha)+b\beta<0
\end{equation}
which is satisfied due to the condition
\[\beta<\frac{b^2-2+2\alpha}{b(2b-1)}.\]
Recall that $\alpha<\frac12+\frac{\beta}{2b}$ from (\ref{alpha}), we thus need to require $\beta< \frac{b(b+1)}{2b^2+b+1}$. Taking $b=1^+$ we have $\beta<\frac12$ and $\alpha<\frac34$.

Now we estimate $\widetilde G_R$. In view of the definition (\ref{m-q1-construct}) of $M_{n+1}$, it follows from the analysis of Subsection \ref{sec-outline}  and Subsection \ref{sec-blocks} that 
\begin{equation}\notag
\begin{split}
\nabla\widetilde G_R=&\ \nabla\widetilde G_n-2\Lambda M_{n+1}\nabla^{\perp}M_{n+1}\\
=&\ \nabla\widetilde G_n-\frac54\lambda_{n+1} \sum_{j=1}^2\xi_j^{\perp}(\xi_j\cdot\nabla)\left(P_{\leq r_{n+1}}a_{j, n+1}\right)^2+\frac52\lambda_{n+1} \nabla^{\perp} M_{n+1}^2\\
&+\nabla J_{NO}+\nabla J_{O1}+\nabla J_{O2}+\nabla J_{O3}+\nabla J_{O4}+\nabla J_{O5}+\nabla J_{O6}
\end{split}
\end{equation}
with $a_{j,n+1}$ replaced by $P_{\leq r_{n+1}}a_{j, n+1}$ and $\lambda_{n+1}$ replaced by $5\lambda_{n+1}$ in the error terms $J_{NO}$, $J_{O1}$, ..., $J_{O6}$.
Denote 
\begin{equation}\notag
\nabla \widetilde G_{R,0}=\nabla\widetilde G_n-\frac54\lambda_{n+1} \sum_{j=1}^2\xi_j^{\perp}(\xi_j\cdot\nabla)\left(P_{\leq r_{n+1}}a_{j, n+1}\right)^2.
\end{equation}
Due to the choice of $a_{j,n+1}$ as in (\ref{a}) to cancel the principal part of $\widetilde G_n$, we have
\begin{equation}\notag
\begin{split}
\nabla \widetilde G_{R,0}=&\ \nabla\widetilde G_n-\frac54\lambda_{n+1} \sum_{j=1}^2\xi_j^{\perp}(\xi\cdot\nabla)P_{\leq 4r_{n+1}}\left(a_{j,n+1}-P_{>r_{n+1}}a_{j, n+1}\right)^2\\
=&\ \frac54\lambda_{n+1} \sum_{j=1}^2\xi_j^{\perp}(\xi_j\cdot\nabla)P_{\leq 4r_{n+1}}\left(2a_{j,n+1}P_{>r_{n+1}}a_{j, n+1}-\left(P_{>r_{n+1}}a_{j, n+1}\right)^2\right)\\
\end{split}
\end{equation}
and hence
\begin{equation}\notag
\begin{split}
\widetilde G_{R,0}=&\ \frac54\lambda_{n+1} \sum_{j=1}^2\Delta^{-1}\nabla\cdot\left(\xi_j^{\perp}(\xi_j\cdot\nabla)P_{\leq 4r_{n+1}}\left(2a_{j,n+1}P_{>r_{n+1}}a_{j, n+1}\right)\right)\\
&-\frac54\lambda_{n+1} \sum_{j=1}^2\Delta^{-1}\nabla\cdot\left(\xi_j^{\perp}(\xi_j\cdot\nabla)P_{\leq 4r_{n+1}}\left(\left(P_{>r_{n+1}}a_{j, n+1}\right)^2\right)\right).
\end{split}
\end{equation}
We further deduce from Lemma \ref{le-aux1} that
\begin{equation}\label{est-gn4}
\begin{split}
\|\widetilde G_{R,0}\|_X\lesssim&\ \lambda_{n+1}\log r_{n+1}\sum_{j=1}^2\|a_{j,n+1}\|_{L^\infty}\|P_{>r_{n+1}}a_{j,n+1}\|_{L^\infty}\\
\lesssim&\ \lambda_{n+1}r_{n+1}^{-2}\log r_{n+1}\sum_{j=1}^2\|a_{j,n+1}\|_{L^\infty}\|\Delta a_{j,n+1}\|_{L^\infty}\\
\lesssim&\ \lambda_{n+1}r_{n+1}^{-2}\log r_{n+1}\sum_{j=1}^2\lambda_n^2\|a_{j,n+1}\|_{L^\infty}^2\\
\lesssim&\ \lambda_{n+1}r_{n+1}^{-2}\lambda_n^2\delta_n\lambda_{n+1}^{-1}\log r_{n+1}\\
\leq&\ \frac{1}{24}\delta_{n+1}
\end{split}
\end{equation}
since $b>1$ and $0<\beta<1$.

We choose
\begin{equation}\notag
\begin{split}
J_{NO}=&-\frac52\lambda_{n+1}\sum_{j=1}^2\Delta^{-1}\nabla\cdot\left(T_{2, 5\lambda_{n+1}\xi_j}[a_{j,n+1}]\xi_j^{\perp} a_{j,n+1}\right)\\
&+\frac12\sum_{j=1}^2\Delta^{-1}\nabla\cdot\left(T_{1, 5\lambda_{n+1}\xi_j}[a_{j,n+1}]\nabla^{\perp} a_{j,n+1}\right).
\end{split}
\end{equation}
Applying Lemma \ref{le-t2} and Lemma \ref{le-aux1} we infer
\begin{equation}\label{est-gn5}
\begin{split}
\|J_{NO}\|_{X}\lesssim &\ \lambda_{n+1}\lambda_{n+1}^{-2+1}\frac{r_{n+1}^2}{\lambda_{n+1}^2}\lambda_{n+1}\log r_{n+1}\|a_{j,n+1}\|_{L^\infty}^2\\
&+\lambda_{n+1}^{-2+1}\frac{r_{n+1}^2}{\lambda_{n+1}}\lambda_{n+1}\log r_{n+1}\|a_{j,n+1}\|_{L^\infty}^2\\
\lesssim&\ \frac{r_{n+1}^2}{\lambda_{n+1}}\frac{\delta_n}{\lambda_{n+1}}\log r_{n+1} \\
\leq&\ \frac{1}{24}\delta_{n+1}
\end{split}
\end{equation}
since $b>1$ and $0<\beta<1$.

For the oscillatory errors we have
\begin{equation}\notag
\begin{split}
J_{O1}=&\ \frac52\lambda_{n+1} \sum_{j=1}^2\Delta^{-1}\nabla\cdot\left( \xi_j^{\perp}(\xi_j\cdot\nabla)a^2_{j,n+1}\cos(10\lambda_{n+1}\xi_j\cdot x) \right)\\
&+\frac52\lambda_{n+1} \sum_{j=1}^2\Delta^{-1}\nabla\cdot\left(T_{2, 5\lambda_{n+1}\xi_j}[a_{j,n+1}]\xi_j^{\perp} a_{j,n+1}\cos(10\lambda_{n+1}\xi_j\cdot x)\right)\\
&+\frac52\lambda_{n+1} \sum_{j=1}^2\Delta^{-1}\nabla\cdot\left( \xi_j^{\perp}(\lambda_{n+1}^{-1} T_{1, 5\lambda_{n+1}\xi_j}[a_{j,n+1}]\nabla^{\perp} a_{j,n+1} \cos(10\lambda_{n+1}\xi_j\cdot x)\right)\\
=:&\ J_{O11}+J_{O12}+J_{O13}.
\end{split}
\end{equation}

\begin{equation}\notag
\|J_{O11}\|_{X}\lesssim \lambda_{n+1}\lambda_{n+1}^{-2+1}\lambda_n\|a_{j,n+1}\|_{L^\infty}^2\lesssim \lambda_n\lambda_{n+1}^{-1}\delta_n.
\end{equation}

Applying Lemma \ref{le-t1} gives
\begin{equation}\notag
\begin{split}
\|T_{1, 5\lambda_{n+1}\xi_j}[a_{j,n+1}]\|_{L^\infty}\lesssim&\ \lambda_{n+1}^{-1}r_{n+1}^2\|a_{j,n+1}\|_{L^\infty}\lesssim \lambda_{n+1}^{-1}r_{n+1}^2\delta_n\lambda_{n+1}^{-1},\\
\|T_{2, 5\lambda_{n+1}\xi_j}[a_{j,n+1}]\|_{L^\infty}\lesssim&\ \lambda_{n+1}^{-2}r_{n+1}^3\|a_{j,n+1}\|_{L^\infty}\lesssim \lambda_{n+1}^{-2}r_{n+1}^3\delta_n\lambda_{n+1}^{-1}.
\end{split}
\end{equation}
Thus we obtain
\begin{equation}\notag
\begin{split}
\|J_{O12}\|_{X}\lesssim&\ \lambda_{n+1}\lambda_{n+1}^{-2+1}\|T_{2, 5\lambda_{n+1}\xi_j}[a_{j,n+1}]\|_{L^\infty}\|a_{j,n+1}\|_{L^\infty}\\
\lesssim&\ \lambda_{n+1}\lambda_{n+1}^{-2+1}\lambda_{n+1}^{-2}r_{n+1}^3\delta_n^2\lambda_{n+1}^{-2}\\
\lesssim&\ \lambda_n\lambda_{n+1}^{-1}\delta_n,
\end{split}
\end{equation}
\begin{equation}\notag
\begin{split}
\|J_{O13}\|_{X}\lesssim&\ \lambda_{n+1}\lambda_{n+1}^{-2+1-1}\|T_{1, 5\lambda_{n+1}\xi_j}[a_{j,n+1}]\|_{L^\infty}\|\nabla^{\perp}a_{j,n+1}\|_{L^\infty}\\
\lesssim&\ \lambda_{n+1}\lambda_{n+1}^{-2+1-1}\lambda_{n+1}^{-1}r_{n+1}^2\lambda_n\delta_n^2\lambda_{n+1}^{-2}\\
\lesssim&\ \lambda_n\lambda_{n+1}^{-1}\delta_n.
\end{split}
\end{equation}
Putting together the estimates above yields
\[\|J_{O1}\|_{X}\lesssim \lambda_n\lambda_{n+1}^{-1}\delta_n. \]
Other oscillatory errors can be estimated analogously and we have 
\begin{equation}\label{est-gn6}
\|J_{O1}\|_{X}+...+ \|J_{O6}\|_{X}\lesssim \lambda_n\lambda_{n+1}^{-1}\delta_n\leq \frac13\delta_{n+1} 
\end{equation}
for $b>1$ and $0<\beta<1$. The estimate (\ref{iter-4}) with $n$ replaced by $n+1$ follows from (\ref{est-gn1})-(\ref{est-gn6}).

Applying (\ref{g1-q1}), (\ref{iter-2}) and the estimates above we also have
\begin{equation}\notag
\begin{split}
\|G_{n+1}\|_{X}\leq&\ \nu\|\Lambda^{\gamma-1} M_{n+1}\|_{X} +\|\widetilde G_N\|_{X}+\|\widetilde G_{R,0}\|_{X}+\|J_{NO}\|_{X}\\
&+\|J_{O1}\|_{X}+...+\|J_{O6}\|_{X}+\|G_n\|_{X}\\
\lesssim &\ \delta_{n+1}+1-\delta_n^{\frac12}\leq 1-\delta_{n+1}^{\frac12}
\end{split}
\end{equation}
for $\lambda_0\gg1$. That is, (\ref{iter-2}) is satisfied with $n$ replaced by $n+1$. Similarly we deduce
\begin{equation}\notag
\begin{split}
\|G_{n+1}\|_{C^s}\leq&\ \nu\|\Lambda^{\gamma-1} M_{n+1}\|_{C^s} +\|\widetilde G_N\|_{C^s}+\|\widetilde G_{R,0}\|_{C^s}+\|J_{NO}\|_{C^s}\\
&+\|J_{O1}\|_{C^s}+...+\|J_{O6}\|_{C^s}+\|G_n\|_{C^s}\\
\lesssim &\ \lambda_{n+1}^{s}\delta_{n+1}+\lambda_{n}^{s}\delta_{n}\\
\lesssim &\ \lambda_{n+1}^{s-\beta}+\lambda_{n}^{s-\beta}\\
\lesssim &\ \lambda_{n+1}^{s}\delta_{n+1}
\end{split}
\end{equation}
for $s\geq \beta$. Thus (\ref{iter-3}) with $n$ replaced by $n+1$ also holds.

In the end, we briefly mention that when $n=q+1$, an analogous analysis as above can be applied to prove the statement of the proposition. The key point is that applying $\eta_{q+1}= \eta_q$, the Nash error term in (\ref{g-q2}) has a better estimate.

\cbdu

\medskip

\subsection{Proof of Theorem \ref{thm}} 
Let $(\theta_0, u_0, f_{1,0})$ and $(\widetilde\theta_0, \widetilde u_0, f_{2,0})$ be two smooth stationary solutions of (\ref{sqg}), with $\theta\neq \widetilde \theta$. Since it is not necessary to require $f_{1,0}\equiv f_{2,0}$, we have the freedom to find two such distinct solutions. Denote
\begin{equation}\notag
\begin{split}
p_0=&\ \frac12(\theta_0+\widetilde\theta_0), \ \ m_0=\frac12(\theta_0-\widetilde\theta_0), \\
\Pi_0=&\ \frac12\Lambda^{-1}(\theta_0+\widetilde\theta_0), \ \ \mu_0= \frac12\Lambda^{-1}(\theta_0-\widetilde\theta_0)
\end{split}
\end{equation}
and let $G_0$ and $\widetilde G_0$ be the stress functions satisfying 
$\Delta G_0=\frac12(f_{1,0}+f_{2,0})$ and $\Delta \widetilde G_0=\frac12(f_{1,0}-f_{2,0})$. One can check $(\Pi_0, \mu_0, G_0, \widetilde G_0)$ satisfies the system (\ref{pm4}). Again we have the flexibility to choose the initial pair of solutions such that (\ref{iter-1})-(\ref{iter-4}) hold for $(\Pi_0, \mu_0, G_0, \widetilde G_0)$. We then apply Proposition \ref{prop} iteratively and obtain a sequence of approximating solutions $\{(\Pi_n, \mu_n, G_n, \widetilde G_n)\}_{n\geq 0}$ of (\ref{pm4}) with $\Pi_n, \mu_n\in C^\alpha$ and $G_n, \widetilde G_n\in C^\beta$. Moreover, $\|\widetilde G_n\|_{X}\leq \delta_n=\lambda_n^{-\beta}$ for $0<\beta<\frac12$. Taking the limit as $n\to\infty$ in the sequence, we obtain a limit solution $(\Pi, \mu, G, 0)$ of (\ref{pm4}) with $\Pi, \mu\in C^\alpha$ and $G\in C^\beta$. We observe $\mu\neq 0$ since the increments are localized around different frequencies. The parameter conditions in Proposition \ref{prop} imply $G\in C^{2\alpha-1}$. Equivalently, $\theta=\Lambda(\Pi+\mu)$ and $\widetilde\theta=\Lambda(\Pi-\mu)$ are two distinct stationary solutions of (\ref{sqg}) with forcing $f=\Delta G$. It completes the proof.

\bigskip

\section*{Acknowledgement}
The authors would like to express their gratitude to Susan Friedlander and Hyunju Kwon for valuable conversations. M. Dai is also grateful for the hospitality of Princeton University ant the Institute for Advanced Study.

%\Endrefs
\end{document}